# An improved convergence theorem for the Newton method under relaxed continuity assumptions

Andrei Dubin

ITEP, 117218, B.Cheremushkinskaya 25, Moscow, Russia

**Abstract**

In the framework of the majorization technique, an improved condition is proposed for the semilocal convergence of the Newton method under the mild assumption that the derivative $F'(x)$ of the involved operator $F(x)$ is continuous. Our starting point is the Argyros' representation of the optimal upper bound for the distance $\|x_{k+1} - x_k\|$ between the adjacent members of the Newton sequence $\{x_k\}$. The major novel element of our proposal is the optimally reconstructed «first integral» approximation $\upsilon_{k+1} = \psi(\upsilon_k)$ to the recurrence relation defining the scalar majorizing sequence $\{\upsilon_k\}$. Compared to the previous results of Argyros, it enables one to obtain a weaker convergence condition that leads to a better bound on the location of the solution of the equation $F(x) = 0$ and allows for a wider choice of initial guesses $x_0$. In the simplest case of the Lipschitz continuous operator $F'(x)$, the new convergence condition improves the famous Kantorovich condition provided that $l_0/l < (6 - 4\sqrt{2})$, where $l_0$ and $l$ stand for the center-Lipschitz and Lipschitz constants respectively.

**Keywords:** Newton method, Banach space, relaxed continuity assumptions, majorization technique, convergence conditions, Kantorovich theorem.

## 1  Introduction

Resolution of many different problems of the fundamental and applied mathematics involves the Newton method generating the corresponding sequence $\{x_k\}$ through the recurrence relation

$$x_{k+1} = x_k - (F'(x_k))^{-1} F(x_k) = T(x_k) \quad , \quad k \geq 0, \tag{1-1}$$

where $F(x) \in Y$ is a continuously differentiable operator which defines a mapping $D_F \subseteq X \to Y$ between an open subset $D_F$ of a Banach space $X$ into a given Banach space $Y$.

The application of this method is built on the so-called Newton-Kantorovich (NK) theorem due to Kantorovich, see [ 1 ] and [ 2 ] for the original proof and further references. The major limitation of this famous theorem is that one has to impose rather stringy smoothness restriction that the operator $F(x)$ is twice differentiable or at least that the Frechet derivative $F'(x)$ of $F(x)$ is Lipschitz continuous. The latter restriction cannot be satisfied in many interesting applied problems. Therefore, the practically important objective is to obtain Kantorovich-type theorems for the Newton method relaxing as much as possible the assumptions about the degree of the continuity of the operator $F'(y)$.

Recently, an interesting new proposal [ 3 ],[ 4 ] in this direction is made presuming only the continuity of $F'(x)$. The employed majorization technique is built on the novel upper bound for the distance $\|x_{k+1} - x_k\|$. The bound is expressed in terms of certain continuity measure depending *only* on the relative distances $\|x_q - x_0\|$ of $x_q$ ($q = k, k-1$) from the initial point $x_0$. The subtle feature which calls for an optimization is that the resulting convergence conditions introduce a larger number of elementary restrictions compared to the NK theorem.

To accomplish such optimization of [ 3 ],[ 4 ], we put forward a new implementation of the majorization technique. It optimally introduces the approximate «first integral» representation [ 2 ] of the recurrence relation defining the majorizing sequence. Under the same computational cost, the proposed technique enables one to improve the convergence condition of [ 4 ] formulating the results under the relaxed continuity assumptions in accordance with the pattern of the NK theorem.

## 2   Improved convergence theorem

### 2.1   «First integral» form of the recurrence relations

To formulate and prove the improved convergence theorem, we derive such implementation of the general majorization condition

$$\|x_{k+1} - x_k\| \leq \upsilon_{k+1} - \upsilon_k \quad , \quad \forall k \geq 0 , \qquad (2\text{-}1)$$

that the scalar majorizing sequence $\{\upsilon_k\}$ is generated by the recurrence relation in the so-called «first integral» form [ 2 ],[ 1 ]:

$$\upsilon_{k+1} = \psi(\upsilon_k) \quad , \quad \forall k \geq 0, \qquad (2\text{-}2)$$

where the continuous non-decreasing function $\psi(\upsilon) \geq 0$ assumes the form

$$\psi(\upsilon) = \eta + 2\bigl(1 - \omega_0(\upsilon)\bigr)^{-1} \cdot \int_0^\upsilon \omega_0(l)\,dl . \qquad (2\text{-}3)$$

Here, the continuous non-decreasing function $\omega_0(r)$ facilitates the affine-invariant upper bound[1]

$$\|(F'(x_0))^{-1} \cdot (F'(x) - F'(x_0))\| \leq \omega_0(r) \quad , \quad \|x - x_0\| \leq r \quad , \quad r \leq R, \qquad (2\text{-}4)$$

where $\bar{B}(x_0, R) \subset D_F$ while $\bar{B}(x_0, R)$ denotes the closed ball of the center $x_0$ and the radius $R$. The (uniform) continuity of $F'(x)$ on $\bar{B}(x_0, R)$ is implemented through the restriction $\omega_0(0) = 0$.

Note that the characteristic feature of the bound ( 2-4 ) is that this continuity measure is *centered* at $x_0$, i.e. depends on the relative distance $\|x - x_0\|$ of $x$ from the *stationary* point identified with $x_0$. Being easier to compute, the continuity measure $\omega_0(\cdot)$ is generically smaller than its conventional non-centered counterpart $\omega(r)$ and, furthermore, $\omega_0(r) \ll \omega(r)$ in a variety of practically interesting cases [ 4 ].

---

[1] The matrix norm is presumed to be submultiplicative and consistent with the corresponding vector norm.

## 2.2 Statement of the theorem

**The set $\Upsilon$.** Consider a continuously differentiable operator $F(x): D_F \subseteq X \to Y$ such that $(F'(x_0))^{-1} \in L(X,Y)$ for some $x_0 \in D_F$. Assume that the condition ( 2-4 ) and the restriction $\eta \geq \|(F'(x_0))^{-1} F(x_0)\| = \|x_1 - x_0\|$ are valid for some constants $\eta > 0$, $R > 0$ and a continuous non-decreasing function $\omega_0(r) \geq 0$ with $\omega_0(0) = 0$. Presume also that ( 2-3 ) defines such continuous non-decreasing function $\psi(\upsilon)$ that the (associated to ( 2-2 )) fixed point equation

$$\upsilon = \psi(\upsilon) \qquad (2\text{-}5)$$

has a minimal solution $\upsilon_*$ satisfying the restrictions $\eta < \upsilon_* \leq R$, $\bar{B}(x_0, R) \subset D_F$ where $\bar{B}(x_0, \upsilon_*)$ denotes the closed ball of the center $x_0$ and the radius $\upsilon_*$. Finally, let the scalar sequence $\{\upsilon_k\}$ be generated by the relation ( 2-2 ) supplemented by the initial conditions $\upsilon_0 = 0$, $\upsilon_1 = \psi(0) = \eta$.

**Theorem.** Under the set $\Upsilon$ of the conditions, the Newton sequence $\{x_k\} \equiv \{T^k x_0\}$ is well-defined, remains in $\bar{B}(x_0, \upsilon_*)$ and converges to a solution $x_*$ of the equation $F(x) = 0$. Moreover, the majorization estimates ( 2-1 ) and

$$\|x_* - x_k\| \leq \upsilon_* - \upsilon_k \quad , \qquad \forall k \geq 0, \qquad (2\text{-}6)$$

are valid where the scalar sequence $\{\upsilon_k\}$, being non-decreasing, converges to the minimal solution $\upsilon_*$ of ( 2-5 ).

In the above stated theorem, the central role is played by the condition that there exists a solution $\upsilon_*$ of the equation ( 2-5 ). To make contact with the standard formulation of the NK theorem, it is sufficient to observe that the NK theorem may be reformulated in terms of the similar requirement. If $2l\eta \leq 1$, the NK fixed point equation $\eta - \upsilon + l\upsilon^2 / 2 = 0$ has a solution $\upsilon_* > \eta$ (for $\eta > 0$) where $l$ denotes the Lipschitz constant in the affine-invariant framework (e.g., see [ 4 ]).

Let us also note without proof that a minor modification of the arguments in [ 4 ] enables one to verify that the considered in the theorem solution $x_*$ is unique in the open ball $B(x_0, \upsilon_*)$.

## 2.3 Proof of the theorem

Assume that there is such a majorizing scalar sequence $\{\upsilon_k\}$ that converges to a limiting point $\upsilon_* = \lim_{k \to \infty} \upsilon_k \leq R$ and facilitates the upper bound ( 2-1 ) once $\{x_k\} \in \bar{B}(x_0, R) \subset D_F$. Then, the standard arguments [ 2 ] verify that the Newton sequence $\{x_k\}$ complies with the Cauchy criterion and, therefore, converges to a limiting point $x_*$ so that the condition ( 2-6 ) is valid and $x_*, x_k \in \bar{B}(x_0, \upsilon_*)$ for $\forall k \geq 0$. In turn, the continuity of $F(x)$ implies that $F(x_*) = 0$.

Once the sequence $\{\upsilon_k\}$ is generated by the recurrence relation ( 2-2 ), the convenient sufficient conditions for the existence of the limiting point $\upsilon_*$ are proposed by Kantorovich [ 1 ].

The proof of Theorem 1 of Section 2 in Chapter XVIII of [ 1 ] includes the verification of the following useful Lemma. Presume that the fixed point equation ( 2-5 ) has a minimal solution $\upsilon_* \in [0, R]$ and the continuous function $\psi(\upsilon) \geq 0$ is non-decreasing when $\upsilon \in [0, R]$. As long as $\{\upsilon_k\}$ is generated by ( 2-2 ) with the considered initial conditions, this sequence is non-decreasing and converges to $\upsilon_*$ so that $\{\upsilon_k\} \in [0, R]$.

Next, let us prove that, under the set of the restrictions $\Upsilon$, the majorization condition ( 2-1 ) is valid where the continuous non-decreasing function $\psi(\upsilon) \geq 0$ is defined on $[0, R]$ by ( 2-3 ). As a starting point of the proof, we use the Rheinboldt majorization technique [ 2 ]. It introduces such a function $\Omega(t, s, r) \geq 0$ that, being continuous and non-decreasing in $t, s, r \geq 0$, *optimally* implements the upper bound $d(T^2 x, Tx) \leq \Omega(d(Tx, x), d(Tx, x_0), d(x, x_0))$ presuming that $\forall x, Tx \in D \subset D_F$. Here, $d(x, y) = \|x - y\|$, $T(x) \equiv Tx$, $T(T(x)) \equiv T^2 x$, $x_0$ denotes the starting point of the Newton sequence $\{x_k\} \in D \subset D_F$ defined by ( 1-1 ) and one identifies $D = \bar{B}(x_0, R)$. In turn, the latter upper bound is converted into the recurrence relation [ 2 ]:

$$u_{k+1} - u_k = \Omega(u_k - u_{k-1}, u_k, u_{k-1}) \quad , \quad \forall k \geq 1. \tag{2-7}$$

Given the initial conditions $u_0 = \upsilon_0 = 0$ and $u_1 = \upsilon_1 = \eta > 0$, it defines such scalar non-decreasing sequence $\{u_k\}$ that facilitates the estimate $\|x_{k+1} - x_k\| \leq u_{k+1} - u_k$ for $\forall k \geq 0$.

As ( 2-7 ) is generically too complex to be analyzed exactly, we propose the following general approach. Given the upper bound $\Omega(\cdot)$, the idea is to *optimally* replace $\Omega(\cdot)$ by such its majorant $\bar{\Omega}(\cdot)$ that admits the so-called «first-integral» representation [ 2 ] (Section 12.5):

$$\Omega(t, t+r, r) \leq \bar{\Omega}(t, t+r, r) = \psi(t+r) - \psi(r) \quad , \quad \forall r \in [0, R] \,, \, \forall t \in [0, R-r], \tag{2-8}$$

where $\psi(r)$ is restricted to be a continuous and non-decreasing function. It is to be compared with the more ambitious suggestion [ 2 ] to try to derive such representation directly for $\Omega(\cdot)$ that is *not* implemented so far for the Newton method under relaxed continuity assumptions.

Given ( 2-8 ) and ( 2-2 ) together with the monotonicity of $\Omega(\cdot)$, it is straightforward to justify by induction that $u_k - u_{k-1} \leq \upsilon_k - \upsilon_{k-1}$ for $\forall k \geq 1$ if $u_0 = \upsilon_0 = 0$ and $u_1 = \upsilon_1 = \psi(0) = \eta \geq 0$. In turn, combining it with the estimate $\|x_{k+1} - x_k\| \leq u_{k+1} - u_k$, one reproduces ( 2-1 ). As for the justification of the inequality $u_k - u_{k-1} \leq \upsilon_k - \upsilon_{k-1}$, note first that it is obviously valid for $k = 1$. Assume that this equality holds true for $1 \leq k \leq q$ which implies that $u_k \leq \upsilon_k$ when $1 \leq k \leq q$. Then, $u_{q+1} - u_q \leq \Omega(\upsilon_q - \upsilon_{q-1}, \upsilon_q, \upsilon_{q-1}) \leq \psi(\upsilon_q) - \psi(\upsilon_{q-1}) = \upsilon_{q+1} - \upsilon_q$ that completes the justification.

It remains to optimally reconstruct the explicit form of the relevant pattern of $\Omega(\cdot)$ and then derive such its majorant $\bar{\Omega}(\cdot)$ that, complying with ( 2-8 ), verifies the expression ( 2-3 ) for $\psi(\upsilon)$. The standard transformations lead to

$$\Omega(t, s, r) = \gamma(s) \cdot \left( \int_r^{r+t} \omega_0(l) dl + \omega_0(r) t \right) \quad , \quad \gamma(s) = (1 - \omega_0(s))^{-1}, \tag{2-9}$$

that, in fact, is the particular case of the general expression obtained in [ 3 ] for a generic Newton-like method (see the equation (34) where one is to identify $a = 0$ and $w_1(r) = w_0(r) = \omega_0(r)$). Thus defined $\Omega(t,s,r)$ is continuous and non-decreasing in $t,s,r$ due to the continuity and monotonicity of $\omega_0(r)$ and the verified below relation $\omega_0(\upsilon_*) \in ]0,1[$. In turn,

$$\Omega(t,s,r) \leq 2\gamma(s)\int_r^{r+t} \omega_0(l)dl \leq 2\gamma(s)\int_0^{r+t} \omega_0(l)dl - 2\gamma(r)\int_0^r \omega_0(l)dl = \bar{\Omega}(t,s,r),  \qquad (\text{2-10})$$

where it is presumed that $r \leq s$ (with $\gamma(r) \leq \gamma(s)$) in accordance with the subsequent identification $r + t = s$. In sum, ( 2-10 ) and ( 2-8 ) lead to ( 2-3 ).

Finally, to prove that thus defined functions $\psi(\upsilon)$ and $\Omega(t,s,r)$ are non-negative, it is sufficient to demonstrate that $\omega_0(\upsilon_*) \in ]0,1[$ which also implies that $\|(F'(x_k))^{-1}F'(x_0)\| < \infty$ for $\forall k \geq 1$ (i.e., the sequence $\{x_k\}$ is well-defined). In turn, this property of $\omega_0(\upsilon_*)$ directly follows from the pattern ( 2-3 ) of ( 2-5 ) and the imposed restriction $\upsilon_* \in ]\eta, R]$.

## 3  Application to the case of the Lipschitz continuous operator $F'(x)$

When $F'(x)$ is Lipschitz continuous, the upper bound ( 2-4 ) is introduced with $\omega_0(r) = l_0 r$ where $l_0$ is the center-Lipschitz constant [ 4 ]. In this case, the implementation ( 2-3 ) of ( 2-5 ) is reduced to the quadratic in $\upsilon$ equation $2l_0\upsilon^2 - (1+l_0\eta)\upsilon + \eta = 0$ with the discriminant $D(l_0\eta) = (1+l_0\eta)^2 - 8l_0\eta = (l_0\eta)^2 - 6l_0\eta + 1 \geq 0$. For the latter equation to possess a solution, the secondary quadratic function $D(\lambda)$ should be non-negative which constrains that $l_0\eta \leq 3 - 2\sqrt{2}$ (i.e., $\eta \leq \eta_{max} = (3-2\sqrt{2})/l_0$) where $\gamma_* = 3 - 2\sqrt{2} \approx 0,171$ is the minimal root of the equation $D(\gamma) = 0$. It is straightforward to demonstrate that the corresponding minimal solution $\upsilon_*$ satisfies the required restriction $\upsilon_* > \eta$ when $\eta > 0$.

The Kantorovich condition [ 1 ] reads $l\eta \leq 0,5$ ($\eta \leq \tilde{\eta}_{max} = (2l)^{-1}$) where $l$ stands for the standard Lipschitz constant [ 4 ] formulated in the affine-invariant way. The above condition $l_0\eta \leq 3 - 2\sqrt{2}$ is weaker if $l_0/l < 2(3-2\sqrt{2})$ when $\eta_{max} = (l_0)^{-1}(3-2\sqrt{2}) > \tilde{\eta}_{max} = (2l)^{-1}$. It is noteworthy that the critical value $6 - 4\sqrt{2} \approx 0,343$ of the ratio $l_0/l$ is fairly moderate because it may often be that $l_0/l \ll 1$ (see [ 4 ]). The price to pay for this advantage is the linear (rather than quadratic as in the Kantorovich case) convergence of the majorizing sequence $\{\upsilon_k\}$.

## 4  Comparison with the Argyros convergence theorem

The proposal of [ 4 ] (Section 2.7) derives certain approximate estimates starting directly from the bound $\Omega(\cdot)$ of ( 2-9 ). The single condition of the existence of the minimal solution $\upsilon_* \in ]\eta, R]$ of the fixed point equation ( 2-5 ) is replaced in [ 4 ] by the two different restrictions. The first one requires the existence of a solution $r_0$ of equation $f(r) = r$ and, in our terms,

$$f(\upsilon) = (1-\omega_0(\upsilon))^{-1} \cdot \left( \int_0^\upsilon \omega_0(l)dl + \omega_0(\upsilon)\upsilon + \eta \right) \geq \psi(\upsilon) \quad , \quad \upsilon \in [0,R]. \tag{4-1}$$

In the derivation of the inequality $f(\upsilon) \geq \psi(\upsilon)$, the monotonicity of $\omega_0(\upsilon)$ is taken into where $\psi(\upsilon)$ is introduced by (2-3). In particular, $f(\upsilon) > \psi(\upsilon)$ for $\forall \upsilon \in ]0,R]$ when $\eta > 0$ and $\omega_0(\upsilon)$ is strictly increasing (with $\omega_0(\upsilon) > 0$ for $\upsilon > 0$). As a result, compared to the first restriction of [4], our condition is weaker leading to a finer localization of the limiting point $x_* \in \bar{B}(x_0, \upsilon_*)$ and a larger maximal admissible value $\eta_{max}$ of the upper bound $\eta \geq \|(F'(x_0))^{-1} F(x_0)\|$ restricting the choice of the initial point $x_0$. As for the second restriction of [4], in our terms it reads $q(r_0) = 2\omega_0(r_0)/(1-\omega_0(r_0)) < 1$ that is *not* generically necessitated by the first restriction.

In particular, given (4-1), the counterpart of our condition $l_0\eta \leq 3 - 2\sqrt{2} \approx 0{,}171$ may be obtained in the form $l_0\eta \leq 0{,}1$ while in [3] it is argued that $l_0\eta \leq (2-\sqrt{3})/2 \approx 0{,}134$ using a slightly different technique. Both these upper bounds are less favorable than our bound.